\documentclass[12pt]{article}
\usepackage{amssymb}
\usepackage{latexsym,bm}
\usepackage{graphicx}
\usepackage{amsmath}
\usepackage{mathrsfs}
\usepackage{amscd,amsthm,xcolor}
\usepackage{tikz}

\usepackage{epstopdf}

\date{}
\newcounter{mathitem}

\DeclareMathOperator{\dist}{dist}
\DeclareMathOperator{\rad}{rad}
\DeclareMathOperator{\ecc}{ecc}
\DeclareMathOperator{\diam}{diam}
	
\begin{document}
\title{Connected graphs with equal independence number and radius}
\author{\hskip -10mm Xingzhi Zhan}
\maketitle
\footnotetext[1]{Department of Mathematics,  Key Laboratory of MEA (Ministry of Education)
 and Shanghai Key Laboratory of PMMP, East China Normal University, Shanghai 200241, China}
\footnotetext[2]{E-mail address: zhan@math.ecnu.edu.cn}
\begin{abstract}
A classical theorem states that the independence number of any connected graph is at least its radius.
In 2007, DeLaVi\~na, Pepper, and Waller proved that every connected graph with equal independence number and radius
is traceable. We first give a new short proof of this theorem and then prove two results for graphs of higher connectivity.
(1) Every 2-connected graph with equal independence number and radius is Hamiltonian. (2) Every 3-connected graph with equal independence
number and radius is Hamiltonian-connected and pancyclic.
Several lemmas used in the proofs are of independent interest.
\end{abstract}

{\bf Keywords.} Hamiltonian graph; traceable graph; pancyclic graph; independence number; radius

{\bf 2020 Mathematics Subject Classification.} 05C45, 05C12, 05C38, 05C69
\vskip 8mm

\section{Introduction and statement of main results}

We consider finite simple graphs and use standard graph-theoretic terminology
as in [16]. The {\it order} of a graph is its number of vertices. We denote by $V(G)$ and $E(G)$ the vertex set and edge set of a graph $G,$ respectively.
For vertices $x,y$ of a connected graph $G$, we denote their distance by
$\dist_G(x,y)$. For a vertex subset $S\subseteq V(G),$ we denote by $G[S]$ the subgraph of $G$
induced by $S.$ Let $K_{s,\, t}$ denote the complete bipartite graph whose partite sets have cardinalities $s$ and $t,$ respectively.

For a connected graph $G$, the eccentricity of a vertex $x$, the radius of
$G$, and the diameter of $G$ are defined, respectively, by
$$
\ecc_G(x)=\max_{y\in V(G)}\dist_G(x,y),
$$
$$
\rad(G)=\min_{x\in V(G)}\ecc_G(x),\qquad
\diam(G)=\max_{x\in V(G)}\ecc_G(x).
$$
The independence number of $G,$ denoted by $\alpha(G),$ is the maximum cardinality of an independent set in $G.$
 A graph is \emph{traceable} if it has a Hamilton path and \emph{Hamiltonian} if it has a Hamilton cycle.

Our starting point is the following classical inequality.

{\bf Theorem 1.} {\it Let $G$ be a connected graph. Then $\alpha(G)\ge \rad(G)$.}

The inequality originated as an early conjecture of Siemion Fajtlowicz's computer program \emph{Graffiti}. It is recorded as
conjecture 0 in \emph{Written on the Wall} [8]. Fajtlowicz and Waller [9] published a direct proof
in 1986. Another independent proof was given by Favaron, Mah\'eo, and Sacl\'e [11]; an earlier contribution of
Favaron appears as a 1986 preprint [10].

A different proof follows from the stronger induced-path theorem of Erd\H{o}s, Saks, and S\'os [6, Theorem 2.2]:
a connected graph of radius $r\ge1$ contains an induced path of order $2r-1.$
Taking alternate vertices of the path gives an independent set of cardinality at least $r$, and hence proves
Theorem 1.

Fajtlowicz's characterization of radius-critical graphs supplies another structural proof [7]. It implies that
every connected graph of radius $r\ge1$ contains an induced \emph{$r$-ciliate}: for some $1\le t\le r$, take a cycle of order $2t$
and attach a path of length $r-t$ at every cycle vertex, with distinct new vertices. When $t=1$, replace the cycle
by a single edge, so the ciliate is a path of order $2r.$ Every $r$-ciliate is bipartite and has $2t(r-t+1)\ge2r$ vertices.
A largest color class is therefore an independent set of cardinality at least $r$. These direct, induced-path,
and radius-critical approaches give several different proofs of Theorem 1.

The equality case in the radius bound was subsequently studied by DeLaVi\~na, Larson, Pepper, and
Waller [4], who gave a structural characterization for $r=1$ and $r\ge5$.

The following theorem was proved by DeLaVi\~na, Pepper, and Waller.

{\bf Theorem 2}  [5]. {\it Let $G$ be a connected graph. If $\alpha(G)=\rad(G)$, then $G$ is traceable.}

Their original proof of Theorem 2 occupies 24 pages and uses 14 lemmas. In this paper we first give a new, short
proof of Theorem 2, and then prove the following two results for graphs with higher connectivity.

{\bf Theorem 3.} {\it Let $G$ be a $2$-connected graph. If $\alpha(G)=\rad(G)$, then $G$ is Hamiltonian.}

For two vertices $x$ and $y$ in a graph, an {\it $(x, y)$-path} is a path with endpoints $x$ and $y.$ A graph $G$ is
{\it Hamiltonian-connected} if for any two distinct vertices $x$ and $y,$ $G$ contains a Hamilton $(x,y)$-path.
A graph $G$ of order $n\ge 3$ is called {\it pancyclic} if for every integer $k$ with $3\le k\le n,$ $G$ contains a cycle of length $k.$

{\bf Theorem 4.} {\it Let $G$ be a $3$-connected graph. If $\alpha(G)=\rad(G)$, then $G$ is Hamiltonian-connected and pancyclic.}

Our proofs share a tree extension lemma, which yields an upper bound on the
independence number of each distance neighborhood. When
$\alpha(G)=\rad(G)$, the bounds for the first two neighborhoods imply that $G$ is claw-free
and exclude a particular configuration of four independent vertices.
For traceability, an induced path of order $2r-1$, where $r=\rad(G)$,
is dominating, and a path-insertion lemma of Brandst\"adt, Dragan,
and K\"ohler [1] completes the proof. For Hamiltonicity, we use
Brousek's induced-subgraph theorem for 2-connected non-Hamiltonian
claw-free graphs [2].

We give two proofs of Hamiltonian-connectedness. The first establishes
distance claw-freeness and applies the theorem of Li and Schelp [13].
The second uses net-freeness when the common parameter is at least four
and a domination bound when it is at most three, together with results
of Shepherd [14] and Vr\'ana, Zhan, and Zhang [15]. Finally, a theorem
of Gould, \L{}uczak, and Pfender [12] yields pancyclicity from the same
local restrictions.

A computer search shows that there are $94694$ connected graphs of order $11$ with independence number equal to radius. Of these, $93853$ are $2$-connected
and $92680$ are $3$-connected.

The remainder of the paper is organized as follows. Section 2 establishes
the common structural lemmas. Section 3 proves Theorems 2 and 3.
Section 4 proves Theorem 4, giving two proofs of Hamiltonian-connectedness
and then a proof of pancyclicity.

\section{A radius bound and independence in neighborhoods}

Let $G$ be a graph. A set $S\subseteq V(G)$ is \emph{dominating} if every vertex outside $S$
has a neighbor in $S$. In this case, we say that {\it $S$ dominates $G.$}
For a vertex $x$ and a nonempty vertex set $S$, put
$$
\dist_G(x,S)=\min_{z\in S}\dist_G(x,z).
$$
We will need the following lemma, which is common to all the arguments.
It yields the general neighborhood inequality below, whose first two cases
provide the local restrictions needed in the subsequent sections.

{\bf Lemma 5.} {\it Let $G$ be a connected graph, let $T$ be a subgraph of $G$ that is a tree,
and let $I\subseteq V(T)$ be independent in $G$. Then
$$
\rad(G)\le\alpha(G)-|I|+\rad(T)+1.
$$
}

{\bf Proof.}
Starting with $H_0=T$ and $S_0=I$, we recursively construct tree
subgraphs $H_i$ of $G$ and sets $S_i\subseteq V(H_i)$ that are
independent in $G$. If $V(H_i)$ does not dominate $G$, connectedness
provides a vertex $q_i$ at distance exactly two from $V(H_i)$.
Indeed, a shortest path from $V(H_i)$ to an undominated vertex
contains such a vertex. Choose a shortest path $a_ip_iq_i$ with
$a_i\in V(H_i)$. The vertices $p_i,q_i$ lie outside $V(H_i)$.
Adjoin these two vertices and only the edges $a_ip_i,p_iq_i$
to $H_i$, obtaining $H_{i+1}$, and put
$$
S_{i+1}=S_i\cup\{q_i\}.
$$
The resulting subgraph $H_{i+1}$ is a tree. Since $q_i$ has no neighbor
in $V(H_i)$, the set $S_{i+1}$ is independent in $G$. Moreover,
$$
\diam(H_{i+1})\le\diam(H_i)+2.
$$

Since $G$ is finite and each step adds two vertices, the process
terminates after some $m\ge0$ steps with $V(H_m)$ dominating $G$.
Thus
$$
|I|+m=|S_m|\le\alpha(G),\qquad
\diam(H_m)\le\diam(T)+2m.
$$
The radius of a tree is half its diameter rounded up, so
$$
\rad(H_m)=\left\lceil\frac{\diam(H_m)}2\right\rceil
\le\left\lceil\frac{\diam(T)}2\right\rceil+m
=\rad(T)+m.
$$
A central vertex of $H_m$ has distance at most $\rad(H_m)+1$ from every
vertex of $G$, because $V(H_m)$ dominates $G$. Consequently,
$$
\begin{aligned}
\rad(G)&\le\rad(H_m)+1\\
&\le\rad(T)+m+1\\
&\le\alpha(G)-|I|+\rad(T)+1.
\end{aligned}
$$
\hfill$\Box$

Taking $T$ to be a single-vertex tree and $I=V(T)$ in Lemma 5 gives
$\rad(G)\le\alpha(G)$, providing another proof of Theorem 1. Thus Lemma 5 is a generalization of Theorem 1.

For a vertex $v$ and an integer $k\ge0$, define
$$
N_k(v)=\{w\in V(G):\dist_G(v,w)=k\},\qquad
B_k(v)=\{w\in V(G):\dist_G(v,w)\le k\}.
$$
Thus $N_1(v)=N(v)$ is the usual neighborhood.
We use the convention $\alpha(G[\emptyset])=0$.

{\bf Lemma 6.} {\it Let $G$ be a connected graph. For every vertex
$v\in V(G)$ and every positive integer $k$,
\begin{equation}\label{eq:neighborhood}
\alpha\bigl(G[N_k(v)]\bigr)
\le\alpha(G)-\rad(G)+\lceil k/2\rceil+1.
\end{equation}}

{\bf Proof.}
If $N_k(v)=\emptyset$, the assertion follows from Theorem 1.
Otherwise choose a maximum independent set $S$ in $G[N_k(v)]$, put
$s=|S|$, and choose a shortest path
$$
p_0p_1\cdots p_k
$$
with $p_0 = v$ and $p_k\in S$. The set
$$
A=\{p_{k-2j}:1\le j\le\lfloor k/2\rfloor\}
$$
is independent because the path is induced. Every vertex in $A$ has
distance at most $k-2$ from $v$, whereas every vertex in $S$ has distance
$k$ from $v$. Since the distances from $v$ of adjacent vertices differ
by at most one, no edge joins $A$ to $S$. Consequently,
$$
I=S\cup A\quad\hbox{is independent},\qquad
|I|=s+\lfloor k/2 \rfloor.
$$
Here $A=\emptyset$ when $k=1$.

Choose a spanning tree $T$ of $G[B_k(v)]$ that preserves distances
from $v$. In fact, such a tree can be found by Dijkstra's algorithm [16, p.97].
Every vertex has distance at most $k$ from $v$ in $T$, so $\rad(T)\le k$.
Since $I\subseteq V(T)$, Lemma 5 gives
$$
\begin{aligned}
\rad(G)&\le\alpha(G)-s-\lfloor k/2 \rfloor+k+1\\
&=\alpha(G)-s+\lceil k/2 \rceil+1.
\end{aligned}
$$
Rearranging proves \eqref{eq:neighborhood}.
\hfill$\Box$

The cases $k=1$ and $k=2$ contain the local information needed below.
The \emph{local independence number} of $G$ is defined by
$$
\alpha_{\ell}(G)=\max_{v\in V(G)}\alpha\bigl(G[N(v)]\bigr).
$$
The bound for $\alpha_{\ell}(G)$ in Corollary 7 below is due to Fajtlowicz;
a proof by DeLaVi\~na and Fajtlowicz appears in [3, Theorem 2].
Here it follows from Lemma 6 together with the corresponding bound
for second neighborhoods.

{\bf Corollary 7.} {\it For every connected graph $G$ and every vertex $v$ of $G$,
$$
\alpha\bigl(G[N_i(v)]\bigr)\le\alpha(G)-\rad(G)+2
\qquad(i=1,2).
$$
In particular,
$$
\alpha_{\ell}(G)\le\alpha(G)-\rad(G)+2.
$$}

{\bf Proof.} Apply Lemma 6 with $k=1,2$, and maximize the first
inequality for $i=1$ over all vertices $v$.
\hfill$\Box$

{\bf Corollary 8.} {\it If a connected graph $G$ contains an independent set of four
vertices $u,v_1,v_2,v_3$ such that
$$
\dist_G(u,v_i)=2\qquad(i=1,2,3),
$$
then $\rad(G)\le\alpha(G)-1$.}

{\bf Proof.}
The vertices $v_1,v_2,v_3$ form an independent set in $G[N_2(u)]$.
Thus Corollary 7 gives
$3\le \alpha(G[N_2(u)])\le\alpha(G)-\rad(G)+2$, as required.
\hfill$\Box$

A \emph{claw} is a graph isomorphic to $K_{1,3}$,
and a graph is \emph{claw-free} if it contains no induced claw.
Let $N_{2,1,1}$ be the graph consisting of a triangle $abca$ and
the three pendant paths $ax_1x_2$, $by$, and $cz$, with no other
edges. A graph is {\it induced $\{X,Y\}$-free} if it contains neither an
induced copy of $X$ nor an induced copy of $Y$.

{\bf Corollary 9.} {\it Every connected graph $G$ with
$\alpha(G)=\rad(G)$ is induced $\{K_{1,3},N_{2,1,1}\}$-free.}

{\bf Proof.}
By Corollary 7, the subgraph induced by the neighborhood of any vertex has independence number
at most two. Hence $G$ is claw-free. In an induced $N_{2,1,1}$,
the set $\{a,x_2,y,z\}$ is independent, and $x_2,y,z$ are all at
distance two from $a$. Such an induced subgraph is excluded by Corollary 8.
\hfill$\Box$

\section{Proofs of Theorems 2 and 3}

We first prove the traceability assertion in Theorem 2 and then the
Hamiltonicity assertion in Theorem 3.

\subsection{Proof of Theorem 2}

A graph is called {\it nontrivial} if it has order at least two.
 A {\it dominating path} in a graph $G$ is a path $P$ whose vertex set is a dominating set of $G$;
 that is, every vertex in $V(G)\setminus V(P)$ has a neighbor on $P$.
The following path-insertion lemma is due to Brandst\"adt, Dragan,
and K\"ohler. We include a proof for completeness.

{\bf Lemma 10} [1, Theorem 3.1]. {\it A claw-free graph with a nontrivial induced dominating path is traceable.}

{\bf Proof.} Let $G$ be a claw-free graph with a nontrivial induced dominating path
$P=p_1p_2\cdots p_m$, with $m\ge2$. Partition the vertices outside $P$ into the sets
$$
X_i=\{x\in V(G)\setminus V(P)|\, \min\{j |\, xp_j\in E(G)\}=i\}
\qquad(1\le i\le m).
$$

For $i\ge2$, the set $X_i$ is a clique. Indeed, two nonadjacent
vertices of $X_i$, together with $p_{i-1}$, would be three independent
neighbors of $p_i$, forming an induced claw.

For $2\le i<m$, every vertex $x\in X_i$ is adjacent to $p_{i+1}$.
Otherwise, $p_{i-1},p_{i+1},x$ would be three independent neighbors
of $p_i$. Consequently, each clique $X_i$, for $2\le i<m$, can be
traversed between $p_i$ and $p_{i+1}$, and $X_m$ can be traversed
after $p_m$ along a path.

It remains to handle $X_1$. All its vertices are adjacent to $p_1$,
so claw-freeness gives $\alpha(G[X_1])\le2$.

We first recall that every connected graph $H$ with $\alpha(H)\le2$
is traceable. Suppose a longest path $Q$ in $H$ is nonspanning.
Neither endpoint has a neighbor outside $Q$. By connectedness, a
nonspanning longest path has at least three vertices. Its endpoints cannot
be adjacent: otherwise $Q$ closes into a cycle, and connectedness
provides an edge from that cycle to an outside vertex, yielding
a longer path. Thus its two endpoints and any outside vertex form
an independent triple, a contradiction.

If $X_1$ is nonempty and $G[X_1]$ is connected, traverse a Hamilton
path of $G[X_1]$, followed by $p_1,p_2$.

If $X_1\ne\emptyset$ and $G[X_1]$ is disconnected, the bound $\alpha(G[X_1])\le2$ implies
that it consists of exactly two cliques, say $A,B$. The vertex $p_2$
is adjacent to every vertex of at least one of them, say $A$.
Otherwise, a nonneighbor of $p_2$ in each clique, together with
$p_2$, would form a claw centered at $p_1$. Start by traversing
all vertices of $B$, then $p_1$, then all vertices of $A$, then $p_2$.
The vertices within each clique may be traversed in any order.
If $X_1$ is empty, simply start with $p_1,p_2$.

Finally, for $i=2,\ldots,m-1$, traverse $X_i$ between $p_i$ and
$p_{i+1}$, omitting any empty set; after $p_m$, traverse $X_m$.
The resulting path contains every vertex of $G$.
\hfill$\Box$

We recall the induced-path theorem mentioned in the introduction.

{\bf Lemma 11} [6, Theorem 2.2]. {\it Every connected graph of radius $r\ge 1$ contains an induced path
of order $2r-1$.}

{\bf Proof of Theorem 2.}
Write $r=\rad(G)=\alpha(G)$. If $r=1$, then $G$ is complete and
hence traceable. Assume $r\ge2$.

By Lemma 11,  $G$ contains an induced path
$$
P=p_1p_2\cdots p_{2r-1}.
$$
Its odd-indexed vertices
$$
\{p_1,p_3,\ldots,p_{2r-1}\}
$$
form an independent set of cardinality $r=\alpha(G)$. Hence $P$ dominates
$G$: a vertex outside $P$ with no neighbor on $P$ could be added
to this independent set, contradicting its maximum cardinality.

By Corollary 9, $G$ is claw-free. Lemma 10 therefore gives a Hamilton path of $G$.
\hfill$\Box$

\subsection{Proof of Theorem 3}

For the Hamiltonian conclusion we use the following theorem of Brousek.

{\bf Lemma 12} [2, Theorem 2A, p.64].  {\it Every 2-connected non-Hamiltonian claw-free graph contains an induced
subgraph obtained from two vertex-disjoint triangles
$$
a_1a_2a_3a_1,\qquad b_1b_2b_3b_1
$$
by connecting each corresponding pair $a_i,b_i$ with a connector, which is either a path
of order at least three or a triangle $a_ib_ic_i a_i$.
Apart from their specified endpoints, the connectors are pairwise
vertex-disjoint and disjoint from the two original triangles.
The resulting induced subgraph has no other edges.}

{\bf Proof of Theorem 3.}
Write $r=\rad(G)$ and $\alpha=\alpha(G)$, so $r=\alpha$. By Corollary 8, the graph $G$ cannot
contain four independent vertices $u,v_1,v_2,v_3$ with
$\dist_G(u,v_i)=2$ for $i=1,2,3$.

By Corollary 9, $G$ is claw-free.

Suppose, for a contradiction, that $G$ is non-Hamiltonian. Apply
Lemma 12 and use its notation for the induced subgraph it
provides.

Assume that one connector, say the first, is a path. Write the path as
$$
a_1p q\cdots b_1,
$$
where $q=b_1$ if this path has exactly two edges. For $i=2,3$, let
$z_i$ be the first internal vertex of connector $i$ starting from
$a_i$; if that connector is a triangle, put $z_i=c_i$. Since the
configuration is induced, the set
$$
\{a_1,q,z_2,z_3\}
$$
is independent. The paths
$$
a_1pq,\qquad a_1a_2z_2,\qquad a_1a_3z_3
$$
show that $q,z_2,z_3$ all have distance exactly two from $a_1$.
This contradicts Corollary 8.

Consequently, all three connectors are triangles. The induced
configuration has nine vertices, with triangles
$$
\{a_1,a_2,a_3\},\qquad
\{b_1,b_2,b_3\},\qquad
\{a_i,b_i,c_i\}\quad(i=1,2,3),
$$
and no additional edges.

Consider the independent set
$$
Q=\{a_1,b_2,c_3\}.
$$
All three pairwise distances between vertices of $Q$ are two, as shown by the shortest paths
$$
a_1a_2b_2,\qquad a_1a_3c_3,\qquad b_2b_3c_3.
$$
If $Q$ does not dominate $G$, connectedness gives a vertex $x$ with
$\dist_G(x,Q)=2$. Choose $u\in Q$ with $\dist_G(u,x)=2$.
Then $Q\cup\{x\}$ is independent, and its three vertices other
than $u$ are all at distance two from $u$, contrary to Corollary 8.
Hence $Q$ dominates $G$.

Since every member of $Q$ is within distance two of $a_1$,
domination gives
$$
r\le\ecc_G(a_1)\le3.
$$
The independence of $Q$ gives $\alpha\ge3$. Together with $r=\alpha$, these
inequalities force
$$
r=\alpha=3.
$$

We now show that $\ecc_G(a_1)\le2$. Otherwise, choose a vertex $w$
with
$$
\dist_G(a_1,w)\ge3.
$$
Then $w$ is nonadjacent to
$$
a_1,a_2,a_3,b_1,c_1.
$$
If $wb_2$ is an edge, the vertex $b_2$ and the three vertices
$w,b_1,a_2$ induce a claw. Similarly, an edge $wb_3$ gives a claw
centered at $b_3$ with leaves $w,b_1,a_3$. Since $G$ is claw-free,
$w$ is therefore also nonadjacent to $b_2,b_3$.

Now
$$
\{w,a_2,b_3,c_1\}
$$
is an independent set of cardinality four, contradicting $\alpha=3$.
Thus $\ecc_G(a_1)\le2$, contrary to $r=3$. This final contradiction
proves that $G$ is Hamiltonian.
\hfill$\Box$

\section{Proof of Theorem 4}

We first give two proofs of the Hamiltonian-connected assertion of
Theorem 4, and then prove its pancyclic assertion. Both proofs of
Hamiltonian-connectedness use the structural results of Section 2.
The first uses distance claw-free graphs; the second uses net-freeness
and domination and does not depend on the distance claw-free results.

\subsection{First proof of Hamiltonian-connectedness}

A graph $G$ is \emph{distance claw-free} if
$\alpha(G[N_i(v)])\le2$ for every vertex $v$ and every positive integer $i$.
This class was introduced by Shepherd [14]. We use the following
criterion of Shepherd, stated in [13], and a theorem of Li and Schelp.

{\bf Lemma 13} [13, Lemma 3, p.187]. {\it A claw-free graph $G$ is distance
claw-free if $\alpha(G[N_2(v)])\le2$ for every vertex $v$.}

{\bf Lemma 14} [13, Theorem 3, p.187]. {\it Every 3-connected distance
claw-free graph is Hamiltonian-connected.}

{\bf Corollary 15.} {\it Every connected graph $G$ satisfying
$\alpha(G)=\rad(G)$ is distance claw-free.}

{\bf Proof.}
By Corollary 9, $G$ is claw-free. Corollary 7 gives
$\alpha(G[N_2(v)])\le2$ for every vertex $v$. Lemma 13 now applies.
\hfill$\Box$

{\bf First proof of the Hamiltonian-connected assertion of Theorem 4.}
By Corollary 15, $G$ is distance claw-free. Since $G$ is 3-connected,
Lemma 14 implies that $G$ is Hamiltonian-connected.
\hfill$\Box$

\subsection{Second proof of Hamiltonian-connectedness}

A \emph{net} consists of a triangle and three additional vertices,
each adjacent to a different vertex of the triangle, with no other
edges. A graph is \emph{net-free} if it contains no induced net.
The exclusion of $N_{2,1,1}$ in Corollary 9 makes the following
structural argument particularly short.

{\bf Lemma 16.} {\it If a connected graph $G$ satisfies
$\alpha(G)=\rad(G)\ge4$, then $G$ is net-free.}

{\bf Proof.}
By Corollary 9, $G$ is claw-free and $N_{2,1,1}$-free.
Suppose that $G$ contains an induced net $H$, with triangle $abca$
and leaves $a',b',c'$ adjacent to $a,b,c$, respectively.
We first show that $V(H)$ dominates $G$.

Otherwise, connectedness gives a path $uyx$ with $u\in V(H)$ and
$\dist_G(x,V(H))=2$. Thus $y,x\notin V(H)$ and $x$ has no neighbor
in $H$. If $y$ has a neighbor among $a,b,c$, relabel so that
$ay\in E(G)$. Delete the edge $bc$ from $H$ and add only the
vertices $y,x$ and the edges $ay,yx$, obtaining a tree $T$.
Every vertex of $T$ is at distance at most two from $a$ in $T$, giving ${\rm rad}(T)\le 2$, and
$\{a',b',c',x\}$ is independent in $G$. Lemma 5 therefore gives
$$
\rad(G)\le\alpha(G)-4+2+1=\alpha(G)-1,
$$
a contradiction.

Consequently, every neighbor of $y$ in $H$ belongs to
$\{a',b',c'\}$. There is at least one such neighbor. There cannot
be two, because these two vertices and $x$ would be three independent
neighbors of $y$, contradicting claw-freeness. Thus $y$ has exactly
one neighbor in $H$, say $a'$. But then $G[V(H)\cup\{y\}]$ is an
induced $N_{2,1,1}$, with pendant paths $aa'y$, $bb'$, and $cc'$.
This contradicts Corollary 9 and proves that $V(H)$ dominates $G$.

Every vertex of $H$ has distance at most two from $a$ in $H$.
Domination now gives $\ecc_G(a)\le3$, contrary to $\rad(G)\ge4$.
Hence $G$ is net-free.
\hfill$\Box$

The condition $\alpha(G)=\rad(G)\ge4$ in Lemma 16 cannot be replaced
by $\alpha(G)=\rad(G)\ge3$: the icosahedral graph has independence
number and radius both equal to three and contains an induced net.
We therefore treat the smaller values separately.

The \emph{domination number} $\gamma(G)$ is the minimum cardinality
of a dominating set. We use the following two results due to Shepherd
and to Vr\'ana, Zhan, and Zhang, respectively.

{\bf Lemma 17} [14, Theorem 2.9]. {\it Every 3-connected claw-free, net-free graph is Hamiltonian-connected.}

{\bf Lemma 18} [15, Theorem 1]. {\it Every 3-connected claw-free graph with domination number at
most three is Hamiltonian-connected.}

{\bf Second proof of the Hamiltonian-connected assertion of Theorem 4.}
Put $r=\alpha(G)=\rad(G)$. By Corollary 9, $G$ is claw-free.
If $r\le3$, a maximum independent set is dominating, since otherwise
an additional vertex could be added to it. Hence
$$
\gamma(G)\le\alpha(G)=r\le3.
$$
Lemma 18 shows that $G$ is Hamiltonian-connected.

If $r\ge4$, Lemma 16 shows that $G$ is net-free. Since $G$ is
3-connected and claw-free, Lemma 17 gives the same conclusion.
\hfill$\Box$

\subsection{Proof of pancyclicity}

The pancyclic assertion of Theorem 4 follows from the same local
restrictions, independently of either Hamiltonian-connectedness proof.
We use the following special case of a theorem of Gould, \L{}uczak, and Pfender.

{\bf Lemma 19} [12, Theorem 1.2, p.186]. {\it Every 3-connected induced $\{K_{1,3},N_{2,1,1}\}$-free graph is pancyclic.}

{\bf Proof of the pancyclic assertion of Theorem 4.}
By Corollary 9, $G$ is induced $\{K_{1,3},N_{2,1,1}\}$-free. Since $G$ is
3-connected, Lemma 19 implies that $G$ is pancyclic.
\hfill$\Box$

\section*{\normalsize Declaration of AI Use}

ChatGPT was used to assist in developing and checking the proofs. The author independently verified all mathematical arguments, 
wrote the paper, and takes full responsibility for its content.

\vskip 5mm
{\bf Acknowledgements.} This research was supported by NSFC Grant No. 12271170 and by the
Science and Technology Commission of Shanghai Municipality under
Grant No. 22DZ2229014.

{\bf Statements and Declarations}

The author has no relevant financial or non-financial interests to disclose.
\end{document}